\documentclass[11pt]{amsart}
\usepackage{geometry}                % See geometry.pdf to learn the layout options. There are lots.
\usepackage[parfill]{parskip}    % Activate to begin paragraphs with an empty line rather than an indent
\usepackage{graphicx}
\usepackage{amssymb}
\usepackage{epstopdf}
\usepackage[all, 2cell]{xy}
\title{On Constructive Axiomatic Method}

\author{Andrei Rodin}

\date{\today}


\thanks{
The work is supported by Russian State Foundations for Humanities, research grant number 13-03-00384. \emph{Acknowledgements}: I thank Sergei Artemov, Ilya Egorychev, Colin McLarty, Marco Panza, Graham Priest, Juha R\"aikk\"a, Vladimir Voevodsky and Noson Yanofsky for very valuable discussions and suggestions.}

\begin{document}
\maketitle

\begin{abstract}
The received notion of axiomatic method stemming from Hilbert  is not fully adequate to the recent successful practice of axiomatizing mathematical theories. The axiomatic architecture of Homotopy type theory (HoTT) does not fit the pattern of formal axiomatic theory in the standard sense of the word.  However this theory falls under a more general and in some respects more traditional notion of axiomatic theory, which I call after Hilbert and Bernays \emph{constructive} and demonstrate using the Classical example of the First Book of Euclid's  \emph{Elements}. I also argue that HoTT is not unique in the respect but represents a wider trend in today's mathematics, which also includes Topos theory and some other developments. On the basis of these modern and ancient examples I claim that the received semantic-oriented formal axiomatic method defended recently by Hintikka is not self-sustained but requires a support of  constructive method. Finally I provide an epistemological argument showing that the constructive axiomatic method is more apt to present scientific theories than the received axiomatic method.  
\end{abstract}

\section{Introduction}
\label{intro}
In a recent paper, which surveys different ways of using axioms in mathematical practice, Dirk Schlimm \cite{Schlimm:2013} provides the following general definition of axiomatic theory and claims that it applies to ``all instances of axiomatics from Euclid to the present'' (op. cit. pp. 43-44): 
 
\begin{quote}
I consider an axiom system to consist of a number of statements (the axioms) and a notion of consequence, which is employed to derive further statements from them; the resulting axiomatic theory includes all statements that are derivable from the axioms and which might also contain terms that are defined from the primitive terms that occur in the axioms (i.e., the primitives of the system). (op. cit., p. 40)
\end{quote}
      
By a ``statement'' the author understands here either a meaningful proposition (the case of \emph{material} aka \emph{contentful} axiom system) or a propositional \emph{form}, which becomes a meaningful proposition through an interpretation (op. cit., p. 44)
\footnote{
Schlimm claims, in full generality, that ``any finite set of axioms can be reduced to a single axiom by conjunction'' (op. cit., p. 67). Since \emph{conjunction} is a logical operation on propositions (and by extension - also on propositional forms) I conclude that in Schlimm's understanding the axioms qualify as propositions or propositional forms. }.   

The above definition is fairly standard. At the first glance it appears that it covers everything that can be reasonably called an axiomatic theory. My aim is to persuade the reader that this is not the case. For this end I shall point to certain theories, which are commonly called axiomatic but which do not fall under the above definition. Then I  shall propose a broader notion of axiomatic theory (and of axiomatic method), which covers these problematic cases. As the reader will see such problematic cases are not marginal and exotic but belong to the mainstream mathematics of different epochs including today's mainstream mathematics. In short, in this paper I argue that the received notion of axiomatic method is biased and propose an improvement.

My  main examples are (i) Euclid's geometrical theory of Book 1 of his \emph{Elements} (see Section 5) and (ii) the recent Homotopy Type Theory (HoTT, see Section 8), which makes part of a new perspective foundations of mathematics called \emph{Univalent Foundations} (UF, see \cite{Voevodsky:2013}). I argue contra Schlimm (see op. cit. p. 50) that if one reads Euclid carefully and avoids any hasty modernization of his text then one is \emph{not} in a position to qualify his geometrical theory as an instance of (contentful) axiomatic theory in Schlimm's sense. Although the overall importance of Euclid's \emph{Elements} for the axiomatic thinking is hardly contestable one may still argue that the particularities of Euclid's axiomatic approach, to which I attract the reader's attention, are archaic  and today may have only historical significance. This is why the case of HoTT-UF, which represents the edge of the current mathematical research,  is more essential for my argument. I shall show that HoTT-UF does not qualify as an axiomatic theory in the received sense of the term either.  Even if the significance of UF for the whole of today's mathematics remains controversial, it is clear that this theory reflects an important part of the mainstream mathematical practice of last decades. I shall argue that the special character of the axiomatic architecture of HoTT represents a wider trend in the 20th century logic and mathematics (Sections 6, 7). The fact that the axiomatic architecture of HoTT also shares certain key features with the old good Euclid's geometry is an additional reason why I believe that the received notion of axiomatic method established in the 20th century is biased.  

The rest of the paper is organized as follows. In the next Section I discuss the received notion of axiomatic method and make explicit some epistemological assumptions, which are usually associated with this method. In Section 3 I overview the continuing discussion on constructive aka \emph{genetic} method started by Hilbert \cite{Hilbert:1900}. In Section 4 I motivate and introduce the notion of emph{constructive}  axiomatic method, which combines the key features of the genetic method and the received axiomatic method. In Section 5  I treat the historical example of Euclid's \emph{Elements}, Book 1 and show that Euclid's axiomatic theory qualifies as constructive in the relevant sense of the term.  In Section 6 I discuss the notion of Curry-Howard correspondence and its role in the Categorical logic. This Section serves me as a necessary historical and theoretical preliminary to Sections 7 and 8 where I treat two modern examples of axiomatic theories: Topos theory and HoTT. In the concluding Section 9 I argue that the constructive axiomatic method has some important epistemological advantages over the received axiomatic method, which make it more apt to prospective applications in physics and other natural sciences. This latter discussion is construed as a critique of the recent account of axiomatic method given by Hintikka \cite{Hintikka:2011}.

\section{The Received Views on the Axiomatic Method: Semantic and Syntactic Approaches}
\label{sec:2}
At an early stage of his life-long work on axiomatic foundation of mathematics (1893-1894) Hilbert describes his notion of well-founded axiomatic theory as follows:

\begin{quote}   
Our theory furnishes only the schema of concepts connected to each other through the unalterable laws of logic. It is left to human reason how it wants to apply this schema to appearance, how it wants to fill it with material. (\cite{Hallett&Majer:2004}, p. 104)
\end{quote}

The two key features of this notion of theory are (i) its schematic character and (ii) its logical grounding. Let me first focus on (ii). Hilbert refers to the ``unalterable laws of logic'' as something definite and somehow known. A weaker assumption which we may attribute to Hilbert and which still allows us to make sense of his words is that the laws of logic are epistemically more reliable than any specific mathematical and scientific knowledge, so these laws should serve as a foundation for all mathematical and scientific theories. 

The assumption of fixity of logic is crucial for understanding the schematic character (i) of Hilbert-style axiomatic theories.  Axioms and theorems of non-interpreted formal theory are \emph{propositional forms}, which admit truth values (and thus turn into full-fledged propositions) through an interpretation. An interpretation of a given propositional form amounts to assigning to terms like ``point'', ``straight line'', etc. certain semantic values, which can be borrowed from another mathematical theory or from some extra-theoretical sources like intuition and experience \cite{Hintikka:2011}. One and the same propositional form (and hence one and the same axiomatic theory), generally, admits multiple interpretations. In order to make sense of the idea of propositional form one needs a further assumption, which is the following. The game of multiple interpretations does not equally concern \emph{all} terms of a given schematic theory. Some terms, namely \emph{logical} terms, have fixed meanings, which (at least in the early versions of Hilbert's axiomatic approach) is supposed to be self-evident. The different treatment of logical and non-logical terms reflects assumption (ii) according to which logical concepts  are epistemically prior. In practice this assumption requires to ``fix logic'' first and then use it as a tool for building various axiomatic theories in a view of their intended models.   

Although Schlimm's definition of axiomatic theory does not commit one to any particular epistemological view on logic, his examples of ``imaginary'', i.e., non-interpreted axioms (\cite{Schlimm:2013}, p. 44) are all of the same sort: they contain some meaningful logical terms and some meaningless (i.e., open to multiple possible interpretations) non-logical terms. This lets me think that the fundamental epistemic assumption (ii) about the logical grounding of axiomatic theories makes proper part of the received view on the axiomatic method as described in Schlimm's paper. 

As it stands Schlimm's definition of axiomatic theory identifies such a theory with a set of ``statements'', i.e., propositions, related by some relation of logical inference aka derivability. Let us now consider the case when the axioms are schematic and formulated in a certain symbolic language with a precise syntax (like in ZF). Then we have at least two possibilities for specifying the wanted relation of logical inference: it can be construed either (i) syntactically as a part of the language (in symbols $A \vdash B$) or (ii) semantically as follows: $B$ follows from $A$  just in case when:  
\begin{quote}
In all interpretations where $A$ is true $B$ is also true. 
\end{quote}
(in symbols $A \models B$;  a more precise definition of semantical consequence is given in Section 9 below). If the two relations coincide extensively on the set of all propositions of some given axiomatic theory $T$ then $T$ is said to be \emph{semantically complete}. 

Whether the syntactic consequence relation (i) or the semantic consequence relation (ii) is epistemically prior is a subject of continuing philosophical controversy. Taking the latter \emph{semantic}  view, one may go further and identify a theory with a class of its models
rather than with any particular axiom system together with all theorems derivable from the given set axioms. (Notice that one and the same class of models can be, generally speaking, axiomatized in many different ways.) In a different paper  Schlimm considers this controversy and  quite rightly stresses the fact that it remains wholly within the scope of the received notion of axiomatic method \cite{Schlimm:2006}. Schlimm makes this claim contra van Fraassen who defends a strong version of semantic view (identifying  a scientific theory with the class of its models) and on this ground denies the relevance of axiomatic method in science altogether (\cite{vFraassen:1980}, p. 65). Unlike Schlimm I don't think that van Fraassen's rejection of the received axiomatic method results from a simple terminological confusion. After all the lack of any significant impact of the axiomatic approach onto the mainstream 20th century physics (unlike the 20th century mathematics) is evident and, in my opinion, needs a better explanation. My tentative explanation is this: the way in which an axiomatic theory (in the received sense of the term) specifies the class of its model is not appropriate (or at least not sufficient) for building physical models. The constructive axiomatic method discussed below in this paper is expected to fix this major shortcoming of the received method. I shall come back to this issue and provide an argument supporting these claims in the concluding Section 9.

\section{Genetic method versus Axiomatic method}
\label{sec:3}
In this Section I overview the continuing discussion on the \emph{constructive} aka \emph{genetic} method of introducing mathematical concepts and building mathematical theories, which often (but not always as we shall shortly see) is opposed to the (received) axiomatic method. I show that the distinction between the methods of theory-building is not quite clear-cut and for this reason needs a revision. 

The current discussion on the \emph{genetic} method goes back to Hilbert's early paper on the concept of number (\cite{Hilbert:1900}, English translation \cite{Hilbert:1996}, where the author explains how the concepts of whole, rational and real numbers can be construed from (or as Hilbert puts this ``engendered [erzeugt] by'') the concept of natural number (provided that this latter concept is taken for granted). This can be done in a number of ways, differences between which are not essential for Hilbert in this context. In particular, real numbers can be construed either as Dedekind cuts or as Cauchy sequences of rational numbers, rational numbers can be construed as pairs of whole numbers numbers, and whole numbers - as pairs of natural numbers. Such procedures Hilbert calls \emph{genetic} and opposes them to \emph{axiomatic} procedures. He remarks that historically the genetic method has been used in arithmetic while the axiomatic method has been used only in geometry. Assuming that 

\begin{quote}
[d]espite the high pedagogic and heuristic value of the genetic method, for the final presentation and the complete logical grounding of our knowledge the axiomatic method deserves the first rank. (\cite{Hilbert:1996}, p. 1093)
\end{quote}

Hilbert develops then an axiomatic theory of arithmetic. 

The distinction between the two methods has been widely discussed in the later literature \cite{Cavailles:1938}, \cite{Demidov:1981}, \cite{Piaget:1956} 
\footnote{
Piaget calls the genetic method \emph{operational}.
}
, \cite{Smirnov:1962},  most recently by Elaine Landry \cite{Landry:2013}.  Landry uses this distinction for making an argument against Feferman's claim according to which axiomatic theories of \emph{categories} (such as EM, ETCS and CCAF, see \cite{McLarty:1992}) cannot serve as an independent foundations of Category theory and of the rest of mathematics because they all allegedly use some prior notions of collection (or class) and operation.  Landry argues back that these notions are used only for  a \emph{genetic} underpinning of Category theory but play no role in the logical axiomatic foundation of this theory \emph{stricto sensu}. Instead of defining a category in terms of \emph{classes} of things called ``objects'' and ``morphisms'' Landry  proposes to think of axiomatic theories of categories as follows:

\begin{quote}
[T]he category-theoretic structuralist simply \emph{begins} by assuming the existence of a system of two sorts of things (namely ``objects'' and ``morphisms'') and then brings these things into relationships with one another by means of certain axioms. (op. cit. p. 44) 
\end{quote}

In a following footnote Landry explains that 

\begin{quote}
What is doing the work here is neither the notion of a system nor the notion of an abstract system, rather it is the Hilbertian idea of a theory as a \emph{schema} for concepts that, themselves, are \emph{implicitly} defined by the axioms. Thus we don't need a ``fixed domain of elements'' [..], we do not need a system as a ``collection'' of elements [..] (op. cit. p. 44)
\footnote{
Interestingly Landry diverges here from a later Hilbert's view according to which a ``fixed'' domain is a hallmark of formal axiomatic method, see the quotes from  \cite{Hilbert&Bernays:1934-1939} In Section 4 below.
 }
\end{quote}

For reasons, which I cannot discuss here I do believe together with Landry that a proper axiomatic construction of Category theory should not \emph{logically} depend on any prior notion of class. However I am not convinced by her argument and do not think that axiomatic theories like EM, ETCS and CCAF (in the received sense of ``axiomatic'') indeed achieve this goal. My objection to Landry is this. (For saving space I shall refer to EM; to ETCS and CCAF my argument applies similarly.)  First of all, I assume (following \cite{Hintikka:2011}) that EM requires a notion of semantic consequence making part of its underlying logic. In order to construe the relation of semantic consequence properly one needs to fix some formal semantical framework (again as a part of the logical machinery). When one uses for this end a Tarski-style model-theoretic semantics it \emph{does} involve some notion of class (collection, universe) of individuals. In this case thinking about categories introduced through EM in terms of classes is \emph{not} just a convenient way of representing categories, which has evident heuristic and pedagogical values, but a proper part of the very axiomatic construction of EM, which belongs to its logical machinery. Thus in this case Landry's argument doesn't go through. Now one may look for a different logical machinery using a different formal semantics, which may allow one to construe the notion of semantic consequence without using the notion of class or similar. However Hilbert's notion of theory as a schema referred to by Landry does not immediately provide such an alternative semantical framework. The only way to get around this objection that I can see is to take a strictly syntactic view on the axiomatic method and thus get rid of the notion of semantic consequence altogether.  Following Hintikka \cite{Hintikka:2011} I believe that this notion is an essential element of the axiomatic method, and so such a price is unacceptable. 
    
The purpose of the above argument is not to settle the axiomatic foundations of Category theory (this topic lies outside of the scope of this paper) but demonstrate that the distinction between the genetic and the axiomatics methods is less clear-cut than it might seem. Whether the class-based definition of category is qualified as axiomatic or genetic depends on one's general ideas about logic and logical semantics, which can be a subject of philosophical controversy. If one takes Tarski-style logical semantics to be a part of logic then the class-based definition of category counts as fairly axiomatic. But Landry's ``category-theoretic structuralist'' sees the same definition as merely genetic.

\section{Constructive Axiomatic Method} 
\label{sec:4} 
In this Section I introduce a generalized method of theory-building, which includes the received axiomatic method as well as the genetic method as special cases. After Hilbert and Bernays I call this generalized method the \emph{constructive axiomatic method}. I also trace here some historical origins of this idea.

\subsection{The Idea} 
The lack of clear-cut distinction between the genetic and the (received) axiomatic method suggests the idea of combining the two method into one. The common ground that I use for combining the two methods is the concept of \emph{rule}. On the one hand,  any \emph{axiomatic} theory comprises rules for handling propositions, which we commonly call \emph{logical} rules (or logical laws)
\footnote{
As it is well-known a deductive system can be formally presented either in Hilbert style (i.e., with many logical axioms and only few rules) or in Gentzen style (with many rules and few axioms). A Gentzen style presentation may have zero axioms. But a Hilbert-style system needs at least one rule (usually \emph{modus ponens}). The notion  of deductive system without rules obviously doesn't make any sense while logical axioms (aka tautologies) are dispensable}.  
On the other hand, any \emph{genetic} theory comprises rules for handling objects of non-propositional types such as arithmetical or geometrical objects. Just as logical rules regulate forming further propositions (theorems) from certain given propositions (including axioms),  genetic rules regulate forming further objects from certain given objects. Now think of a theory, which comprises a set of rules, which apply to propositions as well as to objects of some other types. I shall call such a theory \emph{constructive axiomatic theory} and call the method of building such a theory \emph{constructive axiomatic method}. This general notion of constructive axiomatic theory allows for using different sets of rules. This concerns both logical rules (different axiomatic theories may involve different logical calculi) and non-propositional ``genetic'' rules (different theories may allow for different procedures of object-construction).  The proposed notion of constructive axiomatic theory is \emph{more general} than the received notion of axiomatic theory because a constructive axiomatic theory directly operates with objects of multiple types while an axiomatic theory in the received sense of ``axiomatic'' directly operates only  with propositions. 

One may now remark that propositions typically tell us something \emph{about} certain non-propositional objects  (but may also tell us something about other propositions and about themselves). Thus an axiomatic theory which involves only propositions (axioms and theorems), i.e., an axiomatic theory in the received sense of the word, typically deals with a bunch of non-propositional objects.  The mechanism of \emph{aboutness} which links propositions to non-propositional objects is complicated and not yet fully understood\cite{Yablo:2014}. An essential element of this mechanism is the (semantic) \emph{reference}, which allows for interpreting non-logical terms in (non-interpreted) axioms and theorems as described in Section 2 above. Having this remarkable feature of propositions in mind one may wonder whether there is indeed a need to modify the received notion of axiomatic theory by including rules for non-propositional types on equal footing with logical rules. 

My answer to this question is twofold. The first part of the answer comes from my analysis of mathematical practice: I show that the theory of Book 1 of Euclid's \emph{Elements} (Section 5) and the recent HoTT (Section 8) simply do not qualify as axiomatic theories in the received sense of the term but do qualify as constructive axiomatic theories in the above sense. The second part of the answer derives from an epistemological argument developed in Section 9. I disagree with the view according to which the genetic method of theory building has only certain heuristic and pedagogical significance and argue that it is also needed for construing the concept of semantic consequence properly. On this ground I claim that an axiomatic theory in the received sense of the term always needs to be supported by certain mathematical constructions built genetically.
   
The notion of constructive axiomatic theory provides a better theoretical framework for analyzing this situation than the talk about two different notions of theory (genetic and axiomatic) because it allows one to specify more clearly the interplay between propositional and non-propositional objects and the corresponding rules for operating with these objects. In Section 8 I show that  Martin-L\"of's Type theory (MLTT) can be used as a formal framework for such an analysis.

\subsection{Some Historical Sources} 
The idea of combining the genetic and the (received) axiomatic methods into one method of theory-building is not entirely new. I shall mention here three sources from the earlier literature where this idea is discussed in some form. This list is not supposed to be complete and it doesn't include more recent discussions, which I review in other parts of this paper. 

1) Vladimir Smirnov in his 1962 paper on the genetic method \cite{Smirnov:1962} argues that an adequate modern formalization of genetic theories can be achieved \underline{not} through their axiomatization (by which he understands a presentation of such theories ``in the form of axiomatic calculus where rules of inference [..] control the transition from one proposition to another'')  but rather through a ``direct formalization of recursive techniques, i.e., of algorithmic processes''. For this end, the author argues, one should ``broaden the scope of logic'' by considering ``stratagems of actions'' and ``such forms of thought as prescriptions and systems of prescriptions'' as a properly logical subject-matter (op. cit. p. 275-276). We can see that even if Smirnov keeps the genetic and the axiomatic methods apart he tends to consider them from a common logical point of view. 

2) A similar idea aiming at broadening the scope of logic has been put forward in a different context earlier by Kolmogorov  \cite{Kolmogorov:1932} in the form of ``calculus of problems''. As the author explains in a later note, his 1932 paper

\begin{quote}
 was written with the hope that the logic of solutions of problems would later become a regular part of courses on logic. It was intended to construct a unified logical apparatus dealing with objects of two types - propositions and problems. (\cite{Kolmogorov:1991}, p. 452)
\end{quote}

In Section 6 I shall explain how this Kolmogorov's idea relates to Curry-Howard correspondence and today's Categorical logic including Topos logic. 

3) The relationships between the genetic and the axiomatic methods are described by Hilbert and Bernays in the Introduction to their two-volume \emph{Foundations of Mathematics} \cite{Hilbert&Bernays:1934-1939}, (Eng. tr. \cite{Hilbert&Bernays:2010}). This Introduction first published in 1934 provides a picture, which significantly diverges from what Hilbert says about the same issue in his 1900 paper  \cite{Hilbert:1900} reviewed in the last Section. Here is how Hilbert and Bernays see the relationships between the two methods in 1934: 

\begin{quote}
The term axiomatic will be used partly in a broader and partly in a narrower sense. We will call the development of a theory axiomatic in the broadest sense if the basic notions and presuppositions are stated first, and then the further content of the theory is logically derived with the help of definitions and proofs. In this sense, Euclid provided an axiomatic grounding for geometry, Newton for mechanics, and Clausius for thermodynamics. [..]. [F]or axiomatics in the narrowest sense, the \emph{existential form} comes in as an additional factor. This marks the difference between the \emph{axiomatic method} [in the narrow sense?] and the \emph{constructive} or \emph{genetic} method of grounding a theory. While the constructive method introduces the objects of a theory only as a \emph{genus} of things, an axiomatic theory refers to a fixed system of things [..] given as a whole. Except for the trivial cases where the theory deals only with a finite and fixed set of things, this is an idealizing assumption that properly augments the assumptions formulated in the axioms. (\cite{Hilbert&Bernays:2010}, p. 1a)
\end{quote}

We can see that this time Hilbert includes  the genetic (aka constructive) method as well as the formal (aka schematic) axiomatic method under the scope of the concept of axiomatic method in the ``broadest sense''. What exactly is this new ``broadest'' notion of being axiomatic?  Whether or not Hilbert understands here the genetic method in the same sense as he did this in his 1900 paper?  

The \emph{received} notion of being axiomatic suggests the following reading. Hilbert has changed his earlier idea of genetic method and  in 1934 identified it with the traditional contentful (material) version of the (received) axiomatic method. In this case the ``broadest'' axiomatic method in the above quote is the same as the received axiomatic method. The distinction made by Hilbert and Bernays in the above quote is one between the contentful and the schematic versions of axiomatic method. 
 
The above reading does not seem me convincing. Here is just one reason why. In the same Introduction we find the following interesting remark: 

\begin{quote}
Euclid does not presuppose that points or lines constitute any fixed domain of individuals. Therefore, he does not state any existence axioms either, but only construction postulates. (op. cit. p. 20a)
\end{quote}

The Editors of the English edition of this text (C.-P. Wirth, J. Siekman, M. Gabbay, D. Gabbay) rightly notice here that

\begin{quote}
[u]nless we see construction as a state-changing process, this [Hilbert and Bernays' ] remark is actually not justified: From a static logical point of view, existence sentences are always implied by construction postulates (op. cit. p. 20, footnote 20.3)
\end{quote}

I choose to read the above passage charitably and assume that Hilbert and Bernays \emph{do} see Euclid's geometrical constructions as state-changing processes, i.e., as genuine genetic constructions in the sense of Hilbert's 1900 paper. Now given that (i) Hilbert and Bernays qualify Euclid's theory as axiomatic in the ``broadest'' sense and that (ii) the identification of the ``broadest'' axiomatic method with the received notion of axiomatic method implies that Hilbert and Bernays' remark about Euclid is obviously erroneous, I conclude that (iii) Hilbert and Bernays' ``broadest'' axiomatic method is broader than the received axiomatic method and includes the genetic method in the sense of Hilbert's 1900 paper. Thus if I am right the idea of combined genetic/axiomatic method is already present in \cite{Hilbert&Bernays:1934-1939}. 

 In addition to the general idea and the name I take from Hilbert the following insight.  Recall that in \cite{Hilbert:1900} Hilbert qualifies as genetic the introduction of real numbers through Dedekind cuts (on sets of rational numbers) and through Cauchy sequences of rational numbers. Unless some further precautions are taken (for example, by saying that all Cauchy sequences in question are Turing computable) these ways of introducing new theoretical objects obviously do \emph{not} qualify as constructive in anything like the usual sense of the term - whether this later sense of being constructive is made precise through the concept of Turing computability, through the concept of general recursiveness or in any other way, which determines what does it mean to ``construct'' a mathematical object in terms of certain procedures (like computation) and rules for these procedures. Understandably Hilbert finds the naive genetic ``constructions'' of Dedekind cuts and Cauchy sequences to be epistemically unreliable and proposes instead to build the theory of real numbers axiomatically. Nevertheless I find it useful and interesting to develop (in a modern context) Hilbert's general notion of mathematical construction (in \cite{Hilbert:1900} ), which is not limited by any particular set of ``constructive'' procedures (and rules for such procedures) and may cover, in particular, this naive Dedekind's construction. Notice that the modern received notion of axiomatic theory does not specify by itself any particular logical calculus, which can be used for building such a theory. Similarly the right general notion of being constructive should not, in my view, delimit any particular set of ``constructive'' procedures. This general notion should allow one to specify such procedures differently for different constructive theories.

\section{Euclid: Problems versus Theorems}
\label{sec:5}
In this Section I show that the geometrical theory of Euclid's \emph{Elements}, Book 1, if one read it verbatim and avoids convenient paraphrasing, qualifies as a \emph{constructive} axiomatic theory (in the sense explained in the last Section) but not as an axiomatic theory in the received sense of the word.   

Euclid's theory is based on 5 Axioms (aka Common Notions), 5 Postulates and 23 Definitions. The special historical character of Euclid's Common Notions and Definitions is not relevant to the present discussion, so I leave them apart and focus on Postulates. The first three Postulates are as follows (verbatim after  \cite{Euclid:2011}):

\begin{quote}
[P1:]  to draw a straight-line from any point to any point.

[P2:]  to produce a finite straight-line continuously in a straight-line.

[P3:]  to draw a circle with any center and radius.
\end{quote}

As they stand the three Postulates are not propositions and admit no truth-values. Hence they cannot be axioms in the received sense of the word. In their original form the Postulates cannot be used as premises in logical inferences - if by logical inference one understands an operation that takes some propositions (premises) as its input  and produces some other proposition or propositions (conclusion) as its output
\footnote{
The non-propositional character of Euclid's Postulates has been earlier noticed by A. Szabo, see \cite{Szabo:1978}, p. 230. More recently K. Manders \cite{Manders:1995}, \cite{Manders:2008}  stressed the fact that the Greek ``mathematical demonstrative practice'' as represented by Euclid does not reduce to its propositional content, which is made explicit through ``the notions of mathematical theory and formal proof on which so much recent work in philosophy of mathematics is based'' ( \cite{Manders:2008},  p. 68). This remark shows that Manders' reading of Euclid is compatible with mine. However I don't follow Manders when he describes the non-propositional content of Euclid's proof as a special feature of his ``practice''. Having in mind the notion of constructive axiomatic theory I consider Euclid's non-propositional Postulates as proper elements of his \emph{theory} . I don't deny that the difference between theory and practice makes sense in this context but I disagree with how exactly Manders draws this distinction here.   
}.
In fact Postulates 1-3 are themselves schemes, aka \emph{rules}, of certain basic \emph{operations}, which take some geometrical objects as input and produce some other geometrical objects as output. P1-3 qualify as rules in the same sense of the term in which one usually talks about rules of inference in logic. However P1-3 apply to geometrical  rather than to logical operations. The table below specifies inputs (operands) and  outputs (results) for operations falling under P1-3:

\begin {center}
\begin{tabular}{|l|c|r|}
  \hline
  operation type & input & output \\
  \hline
  P1 & two (different) points & straight segment \\
   \hline
  P2 & straight segment  & (extended) straight segment \\
  \hline
  P3 & straight segment and one of its endpoints & circle  \\ 
  \hline
\end{tabular}
\end {center}

These operations are partly composable in the obvious way: the output of P1-operation is used as input for P2- and P3-operations. This system of operations extended by some further basic operations assumed tacitly
\footnote{
Including the construction of the intersection point of lines in an appropriate position, compare Prop. 1 of Book 1. For further details see \cite{Panza:2012}. The incompleteness of P1-3 has no bearing on my argument.
}
serves Euclid for introducing objects of his theory.  Such an introduction is systematic in the sense that it does not reduce to a simple act of stipulation: it is a procedure, which involves certain elementary operations (including P1-P3) and complex operations obtained through the composition of the elementary operations. As soon as the term \emph{deduction} is understood liberally as a theoretical procedure, which generates some fragments of a given theory from the first principles of this theory, one can say that Euclid's geometrical constructions are deductive. The deductive order of geometrical constructions can be called in this case the \emph{genetic} order. As we shall shortly see in Euclid's theory such a ``constructive deduction'' is tightly related to the more familiar kind of \emph{logical} deduction, which operates with propositions. In what follows we shall see how the same phenomenon arises in Topos theory (Section 7), MLTT and HoTT (Section 8).  

Postulates and Axioms of  are followed by the so-called Propositions. This commonly used title is not found in Euclid's original text where things called by later editors ``Propositions'' are simply numbered but not called by any common name \cite{Euclid:2011}. From Proclus' \emph{Commentary} \cite{Proclus:1873}, (Eng. tr. \cite{Proclus:1970}) written in the 5th century A.D. we learn about the tradition dating back to Euclid's own times (and in fact even to earlier times) of distinguishing between the two sorts of ``Propositions'', namely \emph{Problems} and \emph{Theorems}. Euclid's Theorems by and large are theorems in the modern sense of the word: propositions followed by proofs.  But Problems are something different: they are  \emph{derived rules} for making complex geometrical constructions. Like Postulates Problems admit no truth-values and thus don't qualify as propositions either. Unlike Postulates Problems always require a justification aka \emph{solution}. This is why I call these further rules derived. As we shall now see Euclid's Problems are not solved on the basis of (and hence are not derived from) P1-3 and tacit constructive rules \emph{alone}: the propositional part of Euclid's theory also plays a role in it.        
\footnote{
As an example of Problem consider the initial fragment of Proposition 1, Book 1: ``To construct an equilateral triangle on a given finite straight-line''. It is followed by (i) an appropriate construction and (ii) a proof that the obtained construction is equilateral triangle. For interpretation of this Proposition as a rule see \cite{Panza:2012}, p. 93. Panza says here that this Proposition ``provides a new constructive rule''. I make a stronger claim saying that this Proposition (which is not a proposition in the usual logical sense) \emph{is} this new rule. Saying that I distinguish as usual between a given Problem (= constructive rule) and solution of this Problem, which justifies this rule. 
}.

Indeed what is said so far may suggest that Euclid's theory splits into two separate parts: one consisting of constructive rules derived from Postulates and the other consisting of propositions derived from Axioms. Such a split does not occur for two complementary reasons (which I formulate informally):\\
- (solutions of) non-trivial Problems require (propositional) proofs, which show that the obtained constructions have the required properties; \\
- (proofs of) non-trivial Theorems require constructions, which are conventionally called ``auxiliary''.

This explains why the logical deductive order of Theorems and the genetic order of Problems in Euclid's theory do not exist independently but form a joint deductive order 
\footnote{A more detailed analysis of this deductive structure is found in \cite{Rodin:2014}, ch.2.}.
  
The above observations are fairly straightforward and can hardly remain unnoticed by any attentive reader of Euclid. Why in this case Schlimm \cite{Schlimm:2013} and many other authors tend to think of Euclid's theory as a system of propositions? An obvious explanation is this. All Euclid's Postulates and initial fragments (i.e., bare formulations) of Problems can be easily paraphrased into propositions. This can be done at least in two different ways. The following paraphrases of P1 are self-explanatory: 
\\
\underline{P1m (modal)}: Given two (different) points it is always \emph{possible} to produce a straight segment from one given point to the other given point. \\
\underline{P1e (existential)}: Given two (different) points there \emph{exists} a straight segment having these given points as its endpoint.

P1e instantiates what Hilbert and Bernays call the ``existential form'' used in their formal axiomatic method. The key logical feature of this paraphrase (which it shares with P1m) is this this: it translates  Euclid's non-propositional Postulates and Problems into propositions (in case of P1e - to existential propositions). For further references I shall call a procedure, which aims at replacement of all non-propositional content of a given theory by some suitable propositional content, the \emph{propositional translation} of this theory
\footnote{
In addition to my proposed reading of Postulates 1-3 as formation rules for non-propositional types and the familiar propositional reading  one may also consider reading Postulates 1-3 as specific \emph{functions} from one type of geometrical objects to another. Unlike my proposed reading this latter functional reading doesn't go through for P2, which has no definite outcome. Indeed different operations realized with the same operands may well comply to the same rule (as in the case of rules of chess) while a value of any function  needs to be fixed for any of its argument. A more general historical objection against the functional reading is that the general concept of function is wholly absent in Euclid.   
}. 

There are different ways of thinking about propositional translation of given theory $T$. In particular, this procedure may be conceived of as

(i) a purely linguistic procedure, i.e., a special kind of linguistic paraphrasing, which helps to present the theoretical content of $T$ in the standard propositional form;\\
(ii) a procedure, which separates the theoretical content of $T$ from its non-propositional content, which has no theoretical significance;\\
(iii) a procedure, which transforms $T$ into a new theory $T'$ which, generally, has a different theoretical content. 

Each of (i) and (ii) is sufficient for justifying the identification of Euclid's theory with the result of propositional translation of this theory. I believe that those authors who readily accept such an identification  usually have in their minds some combination of these two assumptions. However in my view both these assumptions are erroneous. (i) is ruled out by the remark that making up a proposition from something, which is not a proposition, is not innocent from a logical point of view. The type difference between propositions and non-propositions matters logically and hence theoretically. Notice that the propositional translation of \emph{Elements}, Book 1, is a matter of its logical reconstruction rather than merely linguistic paraphrasing. It doesn't work as straightforwardly as it does with single Postulates. (ii) depends on a strong epistemological assumption according to which all theoretical knowledge is a propositional knowledge, i.e., the \emph{knowledge-that} \cite{Fantl:2012}. This strong assumption deserves a special discussion, which I postpone until Section 9 where I argue that this assumption is not acceptable either. 

Thus my working understanding of propositional translation will be (iii).  Accordingly I do not identify Euclid's theory of geometry presented in \emph{Elements}, Book 1 with the result of any propositional translation of this theory. Since Euclid's theory in its original form includes a non-propositional content in the form of Problems and Postulates I claim that this theory is an instance of \emph{constructive}  axiomatic theory in the sense of Section 4 but not an instance of a purely propositional axiomatic theory, i.e., of an axiomatic theory in the received sense of the word        
\footnote{
Thanks to Proclus we know that  the idea of propositional translation of Problems into Theorems is very old. From the same source we also know about the contemporary competing idea of translating Theorems into Problems of a special sort (see \cite{Proclus:1970}, p. 63-64). Here is the relevant passage: 
\begin{quote}
Some of the ancients, however, such as the followers of Speusippus and Amphinomus, insisted on calling all propositions ``theorems'', considering ``theorems'' to be a more appropriate designation than ``problems'' for the objects of the theoretical sciences, especially since these sciences deal with eternal things. There is no coming to be among eternals, and hence a problem has no place there, proposing as it does to bring into being or to make something not previously existing - such as to construct an equilateral triangle [..]. Thus it is better, according to them, to say that all these objects exist and that we look on our construction of them not as making, but as understanding them [..] Others, on the contrary, such as the mathematicians of the school of Menaechmus, thought it correct to say that all inquiries are problems but that problems are twofold in character: sometimes their aim is to provide something sought for, and at other times to see, with respect to determinate object, what or of what sort it is, or what quality it has, or what relations it bears to something else.'' (\cite{Proclus:1970}, p. 63-64). 
\end{quote}
Speusippus [408 - 338/9 B.C.] is Plato's nephew who inherited his Academy after Plato's death. Amphinomus is otherwise unknown (except few other references found in the same \emph{Commentary}) Speusippus' contemporary. Menaechmus [380 - 320 B.C.] is mathematician of Plato's circle. The dates of these people suggest that these ideas had been already around for quite a while before Euclid completed his \emph{Elements}  at some point between 300 and 265 B.C.
 .}.
 
\section{Curry-Howard Correspondence and Cartesian Closed Categories}
\label{sec:6}
The present Section is a historical and theoretical preliminary to the following two Sections 7, 8 where I treat modern examples of axiomatic theories. In a nutshell the idea of a Curry-Howard correspondence is given in Kolmogorov's 1932 paper \cite{Kolmogorov:1932} (already discussed in Section 4) where the author establishes that his newly proposed \emph{calculus of problems} has the same algebraic structure as the intuitionistic propositional calculus published earlier in 1930 by Heyting  
\footnote{
Artemov  stresses that the Curry-Howard correspondence doesn't provide a general formal semantics for Kolmogorov's \emph{calculus of problems}, which is known as BHK-semantics (after Brouwer, Heyting and Kolmogorov) \cite{Artemov:2004}. However Curry-Howard serves as a more special constructive \emph{computational}  semantics for the same formal structure (the intuitionsistic calculus). Some historical material concerning the history of Curry-Howard correspondence can be found in  \cite{Cardone&Hindley:2009} and \cite{Seldin:2009} but a more focused historical study on the development of Kolmogorov's logical ideas and of the concept of Curry-Howard correspondence still waits to be written.
}.        
It turns out that this correspondence is extendible onto a large class of symbolic calculi including those, which have been developed independently and apparently for very different purposes.  So there were established a number of correspondences (i.e., of more and less precise isomorphisms) between (proof-related) logical calculi (propositional, first-order or higher-order), on the one hand, and computational calculi (the simply-typed lambda calculus, type systems with dependent types, polymorphic type systems), on the other hand \cite{Sorensen&Urzyczyn:2006}. This led to the so called ``proofs-as-programs and propositions-as-types'' paradigm in logic and computer science where propositions are treated as types of objects (namely of their corresponding proofs) along with other (non-propositional) types of objects belonging the \emph{same} theory
\footnote{
I emphasize that we are talking here about the \emph{same} theory because propositions and proofs of given theory can be also made into objects otherwise, namely, by considering them as objects of a metatheory of the given theory. But this is a different matter.  
}. 
Within this paradigm  building proofs is a special case of building theoretical objects in general. Notice that a constructive axiomatic theory in the sense of Section 4 similarly allows for building systems of propositions (i.e. deriving theorems from axioms) on equal footing with building various non-propositional constructions. 

The Curry-Howard correspondence shows that certain symbolic calculi of different kinds share a common structure. It is natural to ask whether this shared structure can be presented in some invariant way, which would not depend on particularities of syntactic presentations of these symbolic calculi. This problem has been solved in 1963  by means of mathematical Category theory by Lawvere \cite{Lawvere:2004}, who observed that certain categories called Cartesian closed categories (CCC)  ``serve as a common abstraction of type theory and propositional logic''  (\cite{Lawvere:2006}, p. 1). In 1968-1972 this observation has been developed by Lambek into what is now known under the name of Categorical Logic. The three way correspondence between  (i) logical calculi (propositions), (ii) computational calculi (types) and (iii) (objects of) CCC and some other appropriate categories) is called in the literature the \emph{Curry-Howard-Lambek} correspondence \cite{Lambek&Scott:1986}.

Let me now tell a bit more about a relevant part of Lawvere's work that touches upon the issue of axiomatic method directly. Lawvere discovered CCC during his work on an alternative axiomatization of Set theory
\footnote{
Lawvere also had other important motivations for introducing CCC, see \cite{Lawvere:1969a}.
}.
His idea was that the non-logical primitive of standard axiomatic Set theories like ZF, namely the binary membership relation $\in$, was badly chosen. Lawvere suggested to use instead the notion of function and the binary operation (i.e., ternary relation) of composition of functions. The resulting axiomatic theory is known as the Elementary Theory of Category of Sets (ETCS) \cite{Lawvere:1964}, \cite{Lawvere:2005}, \cite{McLarty:1992}.  This proposal may appear radical to one who has habituated oneself to ZF and its likes but as it stands this proposal does not require any modification of the standard formal axiomatic method. The new choice of primitives allowed Lawvere to see that the condition of being CCC makes part of the wanted axiomatic description of the category of sets ($Set$). In this context the further  observation that CCC provides a structural description of Curry-Howard correspondence appears as a bonus. 

Let us see more precisely what happens here. The standard axiomatic approach requires to fix logic first and then use it for sorting out intuitive ideas about sets, numbers, spaces and whatnot. If one encounters then some difficulties of a logical nature one may try to modify the assumed logical principles and try this again.  However in the case of ETCS something different happens: the axiomatization reveals the fact that $Set$ is equipped with its proper \emph{internal} logically-related structure, namely CCC, which is not simply transported from the chosen background logic but emerges as a specific feature of category $Set$. 

Although this fact does not prevent one from construing ETCS as a standard formal axiomatic theory  \cite{McLarty:1992} it suggests a reconsidering the place and the role of logic in this theory. In the next Section we shall see how the similar case of axiomatic Topos theory connects to the idea of constructive axiomatic method.

\section{Topos theory}
\label{sec:7}
The concept of \emph{topos} first emerged in Algebraic Geometry in the circle of Alexandre Grothendieck around 1960 as a far-reaching generalization of the standard concept of topological space and didn't have any special relevance to logic. In his seminal  paper \cite{Lawvere:1970a} Lawvere provided an axiomatic definition of topos called today \emph{elementary} topos
\footnote{
 The title ``elementary'' reflects the fact that Lawvere's definition (unlike Grothendieck's original definition) is expressible in the standard first-order formal language \cite{McLarty:1992}. The concept of elementary topos is more general than that of Grothendieck topos: there are elementary toposes, which are not Grothendieck toposes. A precise axiomatic characterization of Grothendieck toposes is given by \emph{Giraud's axioms}, see  \cite{MacLane&Moerdijk:1992}, p. 575.
}. 
Like ETCS the axiomatic theory of elementary topos does not bring by itself any unusual notion of axiomatic theory. However any systematic exposition of topos theory contains a chapter on the \emph{internal logic} of a topos. In standard textbooks the internal logic is introduced as an extra feature on the top of the basic topos construction \cite{McLarty:1992}, \cite{MacLane&Moerdijk:1992}. As usual it has a syntactic part (Mitchel-B\'enabou language) and a formal semantic part, which interprets the Mitchel-B\'enabou syntax in terms of constructions available in the base topos (Kripke-Joyal semantics). Kripke-Joyal semantics assigns to symbols and syntactic expressions, which have an intuitive logical meaning (logical connectives, quantifiers, truth-values, etc.), explicit semantic values, which otherwise can be called \emph{geometrical} (since the base topos is a generalized space). This is not something wholly unprecedented in the history of the 20th century logic: think, for example, of Tarski's topological semantics for Classical and Intuitionistic propositional logic \cite{Tarski:1956}. However this feature of Kripke-Joyal semantics makes it quite unlike a notion of semantics derived from the idea of interpreting a formal theory by assigning explicit semantic values only to its non-logical elements (compare Section 2 above). 

Internal logic $L_{T}$ can be used for developing further axiomatic theories ``internally'' in given topos $T$. It also provides an additional ``internal'' description of $T$ itself (\cite{McLarty:1992}, Ch. 16). If one looks at Lawvere's \cite{Lawvere:1970a} where the axioms for elementary topos first appeared in the press, one can see that namely the \emph{internal} logical analysis of topos concept allows Lawvere to formulate these axioms: Lawvere observes that the internal logical properties of general topos and the internal logical properties of $Set$ share the CCC property (cartesian closedness) and this observation helps him to obtain the wanted axioms for Topos theory through a modification of ETCS axioms. Although Lawvere didn't use the notion of internal logic of a given category in its explicit form of Mitchel-B\'enabou language, he figured out how standard logical concepts  (propositional connectives, truth-values, modality, etc.) can be  ``encoded'' into this given category: 
  
\begin{quote}
 [A] Grothendieck ``topology'' appears most naturally as a modal operator, of the nature ``it is locally the case that'', the usual logical operators, such as $\forall$, $\exists$, $\Rightarrow$ have natural analogues which apply to families of geometrical objects rather than to propositional functions. [..] [I]n a sense logic is a special case of geometry. ( \cite{Lawvere:1970a}, p. 329)
\end{quote}

I consider this Lawvere's observation as a genuine \emph{logical} analysis of Grothendieck's geometric topos concept  and at the same time observe that it is quite unlike the kind of logical analysis, which led Hilbert to his axiomatic theory of Euclidean geometry \cite{Hilbert:1899}. Unlike Hilbert, Lawvere does not use logic as a ready-made tool (recall Hilbert's ``unalterable laws of logic'') for sorting out certain informal geometrical concepts. Instead Lawvere's analysis reveals a specific logical framework \emph{within} the given geometrical concept. This internal logical framework remained wholly unnoticed by Grothendieck and his collaborators when they first developed the topos concept.        

One may argue that the specific character of Lawvere's logical analysis, which I have just stressed, can matter only when we are talking about the way and the context in which the axioms for the elementary topos have been first obtained but not when we are talking about the axiomatic topos theory in its accomplished form. This argument goes through as far as one uses the received notion of axiomatic theory for distinguishing the final result from the context of its discovery. However this case suggests a reconsidering of this received notion. Indeed Lawvere's dictum  ``logic is a special case of geometry'' is \emph{incompatible} with the construal of Topos theory as axiomatic theory in the received sense of ``axiomatic''. For the standard architecture of axiomatic theories is based on the assumption according to which the scope of logic $L$ used for building some theory $T$ is \emph{broader} than the scope of $T$. By this I mean precisely that $L$ can be used not only for reasoning in $T$ but also for reasoning elsewhere, say, in theory $T'$. Let $T$ be the standard axiomatic Topos theory as in \cite{McLarty:1992}, $L$ be Classical (or Intuitionstic) first-order calculus and $T'$ be ZF (or its Intuitionistic version). I can see no sense in which $L$ can be called as a special case of geometry in this standard setting. Logic $L$ can rather be described here as a common background for a large spectrum of theories, which includes geometrical as well as non-geometrical theories.       

At the same time Lawvere's idea of logic as a part of geometry makes perfect sense if by $T$ one understands a hypothetical \emph{constructive} axiomatic theory of toposes (in the sense of ``constructive'' specified in Section 4 above) and by $L$ - the propositional segment of $T$.  Since $T$ is a geometrical theory (because toposes belong to geometry) its objects are geometrical objects. These geometrical objects are of various types including the type of propositions. Thus propositional objects 
\footnote{As we shall shortly see in MLTT such objects are identified with \emph{proofs} of their corresponding propositions.}
are geometrical objects of sort and logic, by Lawvere's word, is a special case of geometry.  
\footnote{
Lawvere's general ideas about logic are strongly motivated by Hegel's logic and in particular by Hegel's distinction between the \emph{objective} and the \emph{subjective} logic. Hegel's inclusion of  the subjective logic into the objective logic corresponds to Lawvere's inclusion of logic into geometry. For further details see \cite{Rodin:2014}, Section 5.8.  
}.

By the date such a constructive axiomatic theory of toposes doesn't exists in an accomplished form. However I have shown in this Section how Lawvere's work on axiomatization of Topos theory relates to the idea of constructive axiomatic theory
\footnote{
For a discussion on constructive aspects of Category theory and Topos theory see \cite{McLarty:2006}. This discussion is useful but it ignores the aspect of constructivity, which I just have stressed in the main text.   
}.
In the next Section I describe a more recent theory which uses a constructive axiomatic architecture explicitly and thus serves me as a major motivation for this concept (along with historical motivations mentioned in Section 4.2 above).     

\section{Homotopy Type theory and Univalent Foundations}
\label{sec:8}

Homotopy Type theory (HoTT) is a recently emerged field of mathematical research, which has a special relevance to philosophy and logic because it serves as a basis for a new tentative axiomatic foundations of mathematics called the \emph{Univalent Foundations} (UF)\cite{Voevodsky:2011}, \cite{Voevodsky:2013}. In this paper I do not attempt to review HoTT and UF systematically but only describe a special character of its axiomatic architecture.

HoTT emerged through a synthesis of two lines of research, which earlier seemed to be quite unrelated: geometrical Homotopy theory and logical Type theory. The key idea is that of modeling types (including the type of propositions) and terms (including proofs) in Type theory by spaces and their points in Homotopy Theory. Beware that along with basic spaces Homotopy theory also studies \emph{path spaces} where ``points'' are paths in the basic spaces,  spaces of ``paths between paths'' called \emph{homotopies}, spaces of ``paths between paths between paths'' and so on. All these higher-order spaces are also used for interpreting types. 

Like in Topos theory, in HoTT geometry and logic are glued together with appropriate category-theoretic concepts. The central categorical concept used in HoTT is that  of $\omega-groupoid$. In the standard category theory a groupoid is defined as a category where all morphisms between objects are reversible, i.e., are isomorphisms. $\omega$-groupoid is a higher-categorical generalization of this concept (called in this context 1-groupoid) where usual morphisms are equipped by morphisms between morphisms (called 2-morphisms), further morphisms between 2-morphisms (called 3-morphisms) and so on up to the first infinite ordinal $\omega$). An analogy between spaces equipped with paths between points, homotopies between those paths, etc. is straightforward. It allows for mixing the geometrical and the categorical languages and talk interchangeably, e.g.,  about ``spaces of paths'', ``groupoids of paths'' and  ``groupoids'' \emph{simpliciter}. 

The axiomatic HoTT uses Constructive Type theory with depended types due to Martin-L\"of  \cite{Martin-Lof:1984} (MLTT) for turning the tables at this point. The notion of $\omega$-groupoid aka space aka homotopy type is taken as primitive while the notions of proposition, set, (one-dimensional) groupoid, category, etc. are construed as derived notions with MLTT. Types (and, in particular, propositions) and terms (and, in particular, proofs) in MLTT are interpreted, correspondingly, as spaces and points of these spaces. I skip further details. The obtained interpretation of MLTT in the categorically construed Homotopy theory directly translates all constructions in MLTT into geometric constructions. Then one may consider some further axioms such as the  Axiom of Univalence (AU), which gives its name to the Univalent Foundations. In this paper I consider HoTT (with or without AU) only as an axiomatization of modern Homotopy theory. The idea of UF according to which HoTT can be used for developing the rest of mathematics has no bearing on my argument but implies that the special character of HoTT axiomatic architecture may be of general importance for logic and mathematics.  

The constructive character of HoTT is already fully present in MLTT. Unlike ZF and other standard formal axiomatic theories, MLTT is \emph{not} a system of propositions or propositional forms. MLTT construes propositions as a particular type of things among other non-propositional types
\footnote{By propositions I mean here what in \cite{Voevodsky:2013} is called ``mere propositions'':  a type $P$ is a mere proposition if for all $x, y : P$ we have $x = y$ (op.cit. p. 103). In a more general sense each type can be described as a proposition through the Curry-Howard correspondence. 
}
. It comprises rules, which apply to types and their terms in general and to propositions and their proofs in particular. In MLTT, ``the mathematical activity of \emph{proving a theorem} is identified with a special case of the mathematical activity of \emph{constructing an object} - in this case, an inhabitant of a type that represents a proposition'' (\cite{Voevodsky:2013}, p. 17). Since in MLTT the object construction does \emph{not} reduce to proving propositions this theory qualifies as constructive axiomatic theory (in the relevant sense) but not as an axiomatic theory in the received sense. The unusual perspective on logic provided by MLTT has been earlier noticed by Michael Rathjen and described in the following words:  

\begin{quote}
The interrelationship between logical inferences and mathematical constructions connects together logic and mathematics. Logic gets intertwined with mathematical objects and operations, and it appears that its role therein cannot be separated from mathematical constructions. [..] [L]ogical operators can be construed as special cases of more general mathematical operations. ( \cite{Rathjen:2005}, p. 95)\footnote{Notice the analogy with Lawvere's view on logic as a special case of geometry.} 
\end{quote}  

The homotopical interpretation of MLTT (namely, the $\omega$-groupoid interpretation) brings geometry into the picture 
\footnote{
In addition to the homotopical $\omega$-groupoid model HoTT with AU has other ``natural'' models   \cite{Awodey:2014}. Studying and comparing these models from an unified viewpoint largely remains an open research problem. However the $\omega$-groupoid model plays a special role in foundations of mathematics because it supports the idea of reconstructing the world of today's everyday  mathematics using homotopy types as its building blocks. This idea represents what Marquis calls the ``geometrical point of view'' in foundations of mathematics \cite{Marquis:2009}, \cite{Marquis:2013}, which this author opposes to the more familiar Quinean ``logical point of view'' \cite{Quine:1953}. I leave a further discussion on this foundational issue for another occasion.      
} 
and thus makes this case more similar to Euclid's. For example, MLTT validates the following rule (as a special case of a more general rule that I shall not discuss):

(F): Given types $A, B$ to produce type $A \rightarrow B$ of functions with domain $A$ and codomain $B$.

Under the homotopical interpretation this rule becomes

(HF): Given spaces $A, B$ to produce space $A \rightarrow B$ of continuous maps from space $A$ to space $B$,

which has the same form as Euclid's First Postulate allowing for producing a straight line from a given point to another given point
\footnote{
For another analogy between HoTT and Euclid's geometry see \cite{Voevodsky:2013}, p. 56-57. 
}. 
Once again I would like to stress that the non-propositional form of (F) and (HF) just as the non-propositional form of Euclid's Postulates 1-3 is not only a matter of one's favorite style of axiomatic presentation. The case of MLTT/HoTT allows us to see this more clearly than Euclid's case. The type difference in general and the difference between the type of propositions and non-propositional types in particular belongs to the core of MLTT/HoTT, so we cannot  switch here between different types for free. Surely (F) and (HF) can be given a propositional form in a meta-language. However this remark has no bearing on my argument since we are now talking about MLTT and HoTT (i.e., develop an informal metatheory of these theories if one likes) but not discussing the logical form of possible representations of MLTT and HoTT in metatheories of these theories. In MLTT/HoTT propositions has their proper place and rules (F)/(HF) also have their proper place. These places are different. Metatheoretical propositional translations of (F)/(HF) do not belong to MLTT/HoTT and cannot be treated on equal footing with those propositions, which do so. It would be interesting to develop a formal theory of possible metatheoretical propositional translations for rules in MLTT/HoTT and I leave this task for a future research. So far I only claim that the propositional translation of these rules is a meta-theoretical procedure but not just a possible informal interpretation of the same formal structure. 
\footnote{
Compare Section 5 above on propositional translation of Euclid's Postulates 1-3.
}.
Claiming that HoTT is a constructive axiomatic theory I need to pay a particular care for disambiguating the term ``constructive'' properly. MLTT is constructive in a strong sense, which makes this theory \emph{computable}. It is not known to the date whether or not HoTT with AU is constructive in the same strong sense or not; prima facie it is not (\cite{Voevodsky:2013}, p. 11). But in any event HoTT with AU  \emph{is} constructive in the sense of being constructive which I use in this paper. This amounts to saying that  HoTT with AU includes non-propositional types along with the type of propositions. As I have already said HoTT inherits this feature from MLTT. However only in HoTT the non-propositional types are identified as types of certain geometrical objects while in MLTT in its original form such types are presented as a bare type-theoretic abstraction.

\section{What is Constructive Axiomatic Method?} 
\label{sec:9}
In the preceding Sections I motivated, introduced and demonstrated with some examples the concept of constructive axiomatic theory.  In this final Section I would like to evaluate this concept from an epistemological viewpoint.  As my point of departure I shall take Hintikka's recent paper \cite{Hintikka:2011} where the author defends a semantic version of received axiomatic method and analyses its epistemological significance. I shall try to show that this received method of theory-building is not self-sustained (both logically and epistemically) and needs to be supported by a constructive method like one used in HoTT.  

Hintikka:
\begin{quote}
What is crucial in the axiomatic method [..] is that an overview on the axiomatized theory is to capture all and only the relevant structures as so many models of the axioms. (\cite{Hintikka:2011}, p. 72)
\end{quote}

Where these structures come from? Hintikka gives the following answer: 

\begin{quote}
The class of structures that the axioms are calculated to capture can be either given by intuition, freely chosen or else introduced by experience (ib., p. 83)
\end{quote}

One may wonder how a mathematical structure or a structure of some different sort can be ``given'' or ``introduced''  without being construed axiomatically  beforehand. Should we take at this point a Platonistic view according to which mathematical structures in some form exist independently of our axiomatic descriptions of these structures? Hintikka's answer is different. He rejects the notion of intuition understood as an intellectual analogue of sense-perception and describes the relevant notion mathematical intuition as follows:

\begin{quote}
[N]ew logical principles are not dragged [..] by contemplating one's mathematical soul (or is it a navel?) but by active thought-experiment by envisaging different kinds of structures and by seeing how they can be manipulated in imagination. [..] [M]athematical intuition does not correspond on the scientific side to sense-perception, but to experimentation. (ib., p. 78)
\end{quote}

Thus we get the following picture. Building an axiomatic theory is a complicated two-way process; it is a game with Nature (and perhaps also with Society) where raw empirical and intuitive data effect one's axiomatic construction in progress while this construction in its turn effects back one's choice of further data, which become in this way less raw and more structured. Asking where the process starts is the chicken or the egg kind of question. Hintikka's IF logic with its intended game-theoretic semantics provides a precise mathematical model for such games \cite{Pietarinen&Sandu:2004}. 

My concern is about the kind of games that we need to play with Nature for doing science and mathematics. Although yes-no questioning games indeed play an important role in science and perhaps also in mathematics I claim that this is only the top of an iceberg. The main body of this iceberg is filled by mathematical and empirical constructive activities such as designing new experiments. Considering applications of mathematics outside the pure science we may also mention the importance of mathematics for engineering. In order to design a bridge or a particle accelerator one usually plays with certain mathematical \emph{models} of these things, not with formal axioms. 

Since such activities qualify as instances of Hintikka's ``active thought-experiment'' I hardly diverge from Hintikka up to this point. The divergence comes next. I don't grant Hintikka's view according to which the mathematical thought experimentation is, generally, a spontaneous activity, which should be studied by ``empirical psychologists'' rather than logicians and epistemologists (ib. p. 83). I observe that constructive axiomatic theories like Euclid's geometry, Newton's mechanics, Lawvere's axiomatic Topos theory and Voevodsky's HoTT-UF greatly increase one's capacities of mathematical thought experimentation by providing basic elements (like points and straight lines in Euclid) and precise rules (like Euclid's Postulates) for it. I understand that the spontaneity may play a creative role in mathematics and science just like elsewhere but I claim that \emph{typically} a \emph{constructive} axiomatic organization of science and mathematics makes the thought experimentations in these fields richer and more powerful. Compare playing with pieces of wood with playing chess. Such an organized but rather than wholly spontaneous mathematical thought experimentation is typically used for building mathematical models of physical phenomena, designing bridges, particle accelerators etc.  

What is at issue here is how an axiomatic theory controls its semantic contents. Propositional axiomatic theory $T$ (i.e. an axiomatic theory in the standard sense) motivated by certain intended model $M$ controls $M$ through the \emph{truth-evaluation} of its axioms and theorems in $M$. Following Hintikka I shall understand the relation of logical inference in $T$  as the \emph{semantic} consequence relation. Having granted this I claim that such a method of axiomatic control is not effective  because the concept of semantic consequence is highly sensitive to one's basic semantical setting. As long as the semantic consequence is discussed with respect to intuitive structures coming from the air (intuition, free choice, experience) it remains an imprecise intuitive idea itself. I agree with Hintikka that this fact does \emph{not} imply that one has here the choice between appealing to irrational resources and giving up the semantical view on logical consequence altogether (ib., p. 77-78). In order to construe the relation of semantic consequence with a mathematical precision one should fix some \emph{formal semantics}, which allows for doing the truth-evaluation properly (like in the case of Kripke-Joyal semantics for topos logic)
\footnote{
In order to construe a notion of semantic consequence for given formal language $L$ one needs to 

- fix a formal semantics $M$ for that language 

- take a collections $A_{T}$ of well-formed formulas of $L$ (which may express axioms of a given theory $T$ formalized with $L$)

- specify the class $M(A_{T})$ of models of $A_{T}$  by evaluating formulas from $A_{T}$ in $M$

Formula $\phi$ is called a semantic consequence of $A_{T}$, in symbols $A_{T}\models \phi$,  iff $\phi$ is a tautology in $M(A_{T})$.
}. 
For this end one needs to build a basic mathematical model (in the sense of ``model'' used in science rather than Model theory) of the appropriate class of intended structures suggested us by intuition, experience and what not. The yes-no questioning games by themselves cannot bring about such a basic semantic model because, recall, the wanted model is supposed to allow us to make the truth-evaluation properly. Unless such a semantic framework is established one is not in a position to give to yes-no questions definite answers.  The wanted semantic framework cannot be obtained from intuition and experience directly but it can be built  through appropriate mathematical constructive procedures. 

The above analysis suggests that the truth-evaluation is an advanced rather than basic feature of mathematical and other theories
\footnote{
In \emph{that} sense I agree with Hilbert when in \cite{Hilbert:1900} he claims the priority of the propositional axiomatic presentation over the non-propositional genetic presentation (see Section 3 above). However I disagree with the view according to which a final form of propositional axiomatic presentation can and should not depend of any genetic presentation (whether or not this latter view is Hilbert's).         
}. 
Unlike Topos theory, HoTT in its existing form has no resources for doing truth-evaluation internally. But there is a ``general consensus''  that an internal truth-evaluation for HoTT can be construed within a higher-order topos structure in which HoTT would play the role of internal language  (\cite{Voevodsky:2013}, p. 12). 

Mathematics and science not only seek truths and logical relations between those truths but also look for effective methods of doing this and that. These two kinds of knowledge are known in the literature as the \emph{knowledge-that} and the \emph{knowledge-how} \cite{Fantl:2012}. The former kind of knowledge is propositional but the latter is not.  From historical and biological points of view it appears plain that the knowledge-how is a more primitive form of knowledge, which can be reasonably attributed to certain other animals than humans. Propositional knowledge is more advanced and special and it can hardly exist without a developed human-like language. This is why it is plausible that the knowledge-how does not necessarily depend on some knowledge-that. At the same time at least such a developed form of knowledge-that as a propositional axiomatic theory clearly does depend on certain knowledge-how, namely, on the logical knowledge of how to derive theorems from the axioms. In this paper I have shown that this logical knowledge-how in its turn requires knowing how to build certain mathematical construction of non-propositional types.

As far as we are talking about logical and epistemological analysis of \emph{scientific} knowledge there is a tendency to reduce such knowledge to the relevant knowledge-that and separate the associated knowledge-how (perhaps except the purely logical knowledge-how) either into the special domain  \emph{applied} science (and applied mathematics) or into the special domain of \emph{scientific practice}.  In this paper I have shown that such a separation doesn't work for certain mathematical theories because in these cases the two types of knowledge are interlaced at the very atomic level of theoretical reasoning. On the basis of the epistemological considerations presented in this Section I conjecture that the studied cases are not exceptional, so the received notion of axiomatic theory as a bare system of propositions cannot be an adequate theoretical model of an average mathematical theory. The lack of substantial applications of the received form of the axiomatic method in today's mainstream experimental science stressed by Fraassen  back in 1980 \cite{vFraassen:1980}) suggest that the received notion of axiomatic theory does not adequately represent physical and other scientific theories either.

On the positive side I proposed in this paper a notion of constructive axiomatic theory, which can present the knowledge-how (in the form of systems of rules) and the knowledge-that (in the usual propositional form). The epistemological considerations presented above in this Section suggest that this sort of axiomatic theories may be more appropriate in science than the axiomatic theories in the received sense of the term. For recent attempts to apply (what I call) the constructive axiomatic approach in physics see \cite{Schreiber:2014}, \cite{Schreiber&Shulman:2014}. This issue needs a separate discussion, which I leave for another occasion.      
  
\bibliographystyle{plain} 
\bibliography{catax_paper3} 
\end{document}